\documentstyle[12pt]{article}
\catcode`\@=11
\@addtoreset{equation}{section}

\catcode`\@=12
\newtheorem{Theorem}{Theorem}[section]
\newtheorem{Definition}{Definition}[section]
\newtheorem{Proposition}{Proposition}[section]
\newtheorem{Lemma}{Lemma}[section]
\newtheorem{Corollary}{Corollary}[section]
\newtheorem{Note}{Note}[section]
\title{$J-$holomorphic Curves And Periodic
Reeb Orbits
\thanks{Project 19871044 Supported by NSF}}

\author{Renyi Ma\\
Department of Mathematics \\
Tsinghua University \\
Beijing, 100084\\
People's Republic of China\\
rma@math.tsinghua.edu.cn}

\date { }

\begin{document}
\textwidth=125mm
\textheight=185mm
\parindent=8mm
\frenchspacing
\maketitle

\begin{abstract}
 We study the $J-$holomorphic
 curves in the symplectization of the contact manifolds and
 prove that there exists
 at least one periodic Reeb orbits in any closed contact manifold with
 any contact form by using the
 well-known Gromov's nonlinear Fredholm alternative for
 $J-$holomorphic curves. As a corollary, we give a complete
 solution on the well-known Weinstein conjecture.
\end{abstract}

\section{Introduction and results }

A contact structure on a manifold  is a field of a tangent
hyperplanes (contact hyperplanes) that is nondegenerate at any point.
Locally such a field is defined as the field of zeros of a $1-$form
$\lambda $, called a contact form. The nondegeneracy condition
is $d\lambda $ is nondegenerate on the hyperplanes on which $\lambda $
vanishes; equivalently, in $(2n-1)-$space:

$$\lambda \wedge (d\lambda )^{n-1}\neq 0$$

The important example of contact manifold is
the well-known projective cotangent bundles definded as
follows:

Let $N=T^*M$ be the cotangent bundle of the smooth connected compact
manifold $M$. $N$ carries a canonical symplectic structure $\omega
=-d\lambda $ where $\lambda=\sum _{i=1}^{n}y_idx_i$ is the Liouville
form on $N$, see \cite{arg,elg}. Let $P=PT^*M$ be the oriented projective
cotangent bundle of $M$, i.e. $P=\cup _{x\in M}PT^*_xM$. It is well
known that $P$ carries a canonical contact structure induced by the
Liouville form and the projection $\pi :T^*M\mapsto PT^*M$.

Let $(\Sigma ,\lambda )$
be a smooth closed oriented manifold of dimension
$2n-1$ with a contact form $\lambda $.
Associated to $\lambda$ there are two important structures.
First of all the so-called Reed vectorfield $X_\lambda $ defined
by
$$i_X\lambda   \equiv 1, \ \ i_Xd\lambda  \equiv 0$$
and secondly the contact structure $\xi =\xi _{\lambda }
\mapsto \Sigma $ given by
$$\xi _{\lambda }=\ker (\lambda )\subset T \Sigma $$
by a result of Gray, \cite{gra} , the contact structure is very
stable. In fact, if $(\lambda  _t  )_{t\in  [0,1]}$ is a smooth
arc of contact forms inducing the arc of contact structures
$(\xi _t)_{t\in [0,1]}$, there exists a smooth arc
$(\psi _t)_{t\in [0,1]}$
of diffeomorphisms with $\psi _0=Id$, such that
\begin{equation}
T\Psi _t(\xi _0)=\xi _t
\end{equation}
here it is important that $\Sigma $ is compact. From (1.1) and
the fact that $\Psi _0=Id$ it follows immediately that there
exists a smooth family of maps $[0,1]\times \Sigma \mapsto
(0,\infty ):(t, m)\to f_t(m)$ such that
\begin{equation}
\Psi ^*_t\lambda _t=f_t\lambda _0
\end{equation}
In contrast to the contact structure the dynamics of the
Reeb vectorfield changes drastically under small perturbation
and in general the flows associated to $X_t$ and $X_s$ for
$t\neq s$ will not be conjugated, see\cite{arg,elg,ehs,ho}.

\vskip 3pt

  Let $M$ be a Riemann manifold with Riemann metric, then it is well known
that there exists a canonical contact structure in the unit
sphere of its tangent bundle and the motion of geodesic line lifts
to a geodesic flow on the unit sphere bundles. Therefore
the closed orbit of geodesic flow or Reeb flow
on the sphere bundle projects to a closed geodesics in
the Riemann manifolds, conversely the closed geodesic orbit lifts
to a closed Reeb orbit. The classical work of Ljusternik and Fet
states that every simply connected Riemannian
manifold has at least one closed geodesics, this with the
Cartan and Hadamard's results on non-simply closed Riemann
manifold implies that any closed
Riemann manifolds has a closed geodesics, i.e.,
the sphere bundle of a closed Riemann manifold with
standard contact form carries at least one closed Reeb orbits
which is a lift of closed geodesics of base manifold.
Its proof depends on the classical minimax principle
of Ljusternik and Schnirelman or minimalization
of Hadamard and Cartan,\cite{kl}. In sympletic geometry, Gromov \cite{gro}
introduces the global methods to proves the existences of
symplectic fixed points or periodic orbits which depends on the
nonlinear Fredholm alternative of $J-$holomorphic curves
in the symplectic manifolds. In this paper we use the
$J-$holomorphic curve's method to prove
\begin{Theorem}
Every closed contact manifold $\Sigma $ with contact form $\lambda $
carries at least one closed orbit.
\end{Theorem}
Theorem1.1 was conjectured by Weinstein in \cite{we2} under the assumption
$H^1(\Sigma )=0$. Note that Viterbo \cite{vi} first proved the
above result for any contact manifolds $\Sigma $ of
restricted type in $R^{2n}$ after
Rabinowitz \cite{ra} and Weinstein \cite{we1,we2}.
After Viterbo's work many results
were obtained in \cite{fhv,ho,hz,hv,ma1,ma2} etc
by using variational method or Gromov's
nonlinear Fredholm alternative, see survey paper
\cite{el}. Through $J-$holomorphic curves, especially,
Hofer in \cite{ho} proved the following first striking results.
\begin{Corollary}
(Hofer)The three dimentional sphere with
any contact form carries at least one closed Reeb orbit.
\end{Corollary}

\begin{Theorem}
(Ljusternik-Fet) Every simply connected closed Riemannian
manifold has at least one closed geodesics.
\end{Theorem}
Therefore we get a new proof on the well-known Ljusternik-Fet
Theorem without using the classical minimax principle,
an alternative proof can be found in \cite{ma4}.

\begin{Theorem}
(Cartan-Hadamard) Every non-simply closed
Riemannian manifold has at least one closed
geodesics.
\end{Theorem}
Our method can not conclude that the geodesics
is minimal.

\vskip 3pt

{\bf Sketch of proofs}: We work in the framework
as in \cite{gro,ma3}. In Section 2, we study the
linear Cauchy-Riemann operator and sketch some basic
properties. In section 3, first we construct a
Lagrangian submanifold $W$ under the assumption that
there does not exists closed Reeb orbit in
$(\Sigma ,\lambda )$; second, we study
the space ${\cal D}(V,W)$
of contractible disks in manifold $V$ with boundary
in Lagrangian submanifold $W$ and construct a Fredholm
section of tangent bundle of  ${\cal  D }(V,W)$.
In section 4, following \cite{gro,ma3}, we prove that the Fredholm
section is not proper by using a
special anti-holomorphic section as in \cite{gro,ma3}.  In section 5-6,
we use a geometric argument to deduce the
boundary $C^0-$estimates on $W$. In the final section, we
use nonlinear Fredholm trick in \cite{gro,ma3}
to
complete our proof.

\vskip 3pt

Since the proofs in this paper are very difficult, we suggest the reader
first read the Gromov's paper\cite{gro},
Audin and Lafondaine's book\cite{aul}, and
Hummel's book\cite{hu}.

\begin{Note}
The related problem with Weinstein conjecture(see\cite{we2}) is
Arnold chord conjecture(see\cite{ar}) which was discussed in
\cite{ar,gi} and finally solved in \cite{ma3}.
The generalized Arnold conjecture corresponding
to Theorem 1.1 was also solved in similar method of this paper.
These results was reported in the Second
International Conference on Nonlinear Analysis, 14-19 June
1999, Tianjin, China; First $3\times 3$ Canada-China
Math Congress, Tsinghua University, Beijing, August 23-28,1999;
Differential Geometry Seminar in Nankai Institute of Mathematics, Oct. 24-31, 2000;
Symplectic Geometry Seminar In Nankai Institute of Math., Dec. 28-31, 2000;
International conference on Symplectic geometry in Sichun Uni., June 24-July, 2001.
Some technique part of proofs was carried in ICTP from August to October, 2001.
The auther is deeply grateful to thank for the
all inviters, especially to Professor Y. M. Long.
\end{Note}

\section{Linear Fredholm Theory}

For $100<k<\infty $ consider the Hilbert space
$V_k$ consisting of all maps $u\in H^{k,2}(D, C\times C^n)$,
such that $u(z)\in \{izR\}\times R^n\subset C\times C^n$
for almost all $z\in
\partial D$. $L_{k-1}$ denotes the usual Sobolev
space
$H_{k-1}(D, C\times C^n)$. We define an operator
$\bar \partial :V_k\mapsto L_{k-1}$ by
\begin{equation}
\bar \partial u=u_s+iu_t
\end{equation}
where the coordinates on $D$ are
$(s,t)=s+it$, $D=\{ z||z|\leq 1\}  $.
The following result is well known(see\cite[p96,Th3.3.1]{wen}).
\begin{Proposition}
$\bar \partial :V_k\mapsto L_{k-1}$ is a surjective real
linear Fredholm operator of index $n+3$. The kernel
consists of $(a_0+isz-\bar a_0z^2,s_1,...,s_n)$,
$a_0\in C$, $s,s_1,...,s_n$.
\end{Proposition}
Let $(C^n, \sigma =-Im(\cdot ,\cdot ))$ be the standard
symplectic space. We consider a real $n-$dimensional plane
$R^n\subset C^n$. It is called Lagrangian if the
skew-scalar product of any two vectors of $R^n$ equals zero.
For example, the plane $\{(p,q)|p=0\}$ and $\{(p,q)|q=0\}$ are
two transversal Lagrangian
subspaces. The manifold of all (nonoriented) Lagrangian subspaces of
$R^{2n}$ is called the Lagrangian-Grassmanian $\Lambda (n)$.
One can prove that the fundamental group of
$\Lambda (n)$ is free cyclic, i.e.
$\pi _1(\Lambda (n))=Z$. Next assume
$(\Gamma (z))_{z\in \partial D}$ is a smooth map
associating to a point $z\in \partial D$ a Lagrangian
subspace $\Gamma (z)$ of $C^n$, i.e.
$(\Gamma (z))_{z\in \partial D}$ defines a smooth curve
$\alpha $ in the Lagrangian-Grassmanian manifold $\Lambda (n)$.
Since $\pi _1(\Lambda (n))=Z$, one have
$[\alpha ]=ke$, we call integer $k$ the Maslov index
of the closed curve $\alpha $ and denote it by $m(\Gamma )$
(see\cite[p116-118]{arg}). Note that the Maslov index of the closed
curve $\alpha $ is just the two times of the rotation
numbers(see\cite[p116-118]{arg} or \cite[p96,Th3.3.1]{wen}).

Now let $z:S^1\mapsto R\times R^n\subset C\times C^n$ be a smooth curve.
Then
it defines a constant loop $\alpha $ in Lagrangian-Grassmanian
manifold $\Lambda (n+1)$. This loop defines
the Maslov index $m(\alpha )$ of the map
$z$ which is easily seen to be zero.

  Now Let $(V,\omega )$ be a symplectic manifold,
$W\subset V$ a closed Lagrangian submanifold. Let
$(\bar V,\bar \omega )=(D\times V, \omega _0+\omega )$ and
$\bar W=\partial D\times W$. Let
$\bar u=(id,u):(D,\partial D)\to (D\times V,\partial D\times W)$
be a smooth map homotopic to the map
$(id,u_0)$, here $u_0:(D, \partial D)\to p\in W\subset V$.
Then $\bar u^*TV$ is a symplectic vector bundle on $D$ and
$(\bar u|_{\partial D})^*T\bar W$ be a Lagrangian subbundle in
$\bar u^*T\bar V|_{\partial D}$. Since $\bar u:(D,\partial D)\to
(\bar V,\bar W)$ is homotopic to
$\bar u_0$, here$u_0(z)=(z,p)$, i.e., there exists
a homotopy $h:[0,1]\times (D,\partial D)\to (\bar V,\bar W)$ such
that $h(0,z)=(z,p),h(1,z)=\bar u(z)$, we can take
a trivialization of the symplectic
vector bundle $h^*T\bar V$ on $[0,1]\times (D,\partial D)$ as
$$\Phi (h^*T\bar V)=[0,1]\times D\times C\times C^n$$
and
$$\Phi ((h|_{[0,1]\times \partial D})^*T\bar W)\subset [0,1]\times
S^1\times C\times C^n$$
Let
$$\pi _2: [0,1]\times D\times C\times C^n\to C\times C^n$$
then
$$\tilde h: (s,z)\in [0,1]\times S^1
\to \pi _2\Phi (h|_{[0,1]\times \partial D})^*T\bar W|(s,z)\in \Lambda (n+1).$$

\begin{Lemma}
Let $\bar u: (D,\partial D) \rightarrow (\bar V,\bar W)$
be a $C^k-$map $(k\geq 1)$
as above. Then,
$$m(\tilde u)=2.$$
\end{Lemma}
Proof.  Since $\bar u$ is homotopic to
$\bar u_0$ in $\bar V$ relative to $\bar W$,
by the above argument we have
a homotopy $\Phi _s$ of trivializations
such that
$$\Phi _s(\bar u^*TV)=D\times C\times C^n$$
and
$$\Phi _s((\bar u|_{\partial D})^*T\bar W)\subset S^1\times C\times C^n$$
Moreover
$$\Phi _0(\bar u|_{\partial D})^*T\bar W=S^1\times {izR}\times R^n$$
So, the homotopy induces
a homotopy $\tilde h$ in Lagrangian-Grassmanian
manifold. Note that $m(\tilde h(0, \cdot ))=0$.
By the homotopy invariance of Maslov index,
we know that $m(\tilde  u|_{\partial D})=2$.

\vskip 5pt

   Consider the partial differential equation
\begin{eqnarray}
&&\bar \partial \bar u+A(z)\bar u=0  \ on \ D  \cr
&&\bar u(z)\in \Gamma (z) ({izR}\times R^n)\ for \ z\in \partial D \cr
&&\Gamma (z)\in GL(2(n+1),R)\cap Sp(2(n+1))\cr
&&m(\Gamma )=2 \ \ \ \ \ \ \ \
\end{eqnarray}

For $100<k<\infty $ consider the Banach space $\bar V_k $
consisting of all maps $u\in H^{k,2}(D, C^n)$ such
that  $u(z)\in \Gamma (z)$ for almost all $z\in
\partial D$. Let $L_{k-1}$ the usual Sobolev space $H_{k-1}(D,C\times C^n)$

\begin{Proposition}
$\bar \partial : \bar V_k \rightarrow L_{k-1}$
is a real linear Fredholm operator of index n+3.
\end{Proposition}

\section{Nonlinear Fredholm Theory}

\subsection{Constructions of Lagrangian submanifolds}

Let $(\Sigma ,\lambda )$ be a contact manifolds with contact form
$\lambda $ and $X$ its Reeb vector field, then
$X$ integrates to a Reeb flow $\eta _t$ for $t\in R^1$.
Consider the form $d(e^a\lambda )$
at the point $(a,\sigma )$ on the manifold
$(R\times \Sigma )$, then one can check
that $d(e^a\lambda )$ is a symplectic
form on $R\times \Sigma $. Moreover
One can check that
\begin{eqnarray}
&&i_X(e^a\lambda )=e^a \\
&&i_X(d(e^a\lambda ))=-de^a
\end{eqnarray}
So, the symplectization of Reeb vector field $X$ is the
Hamilton vector field of $e^a$ with
respect to the symplectic form $d(e^a\lambda )$.
Therefore the Reeb flow lifts to the Hamilton flow
$h_s$ on $R\times \Sigma $(see\cite{arg,elg,ehs}).

Let
$$(V',\omega ')=((R\times \Sigma )\times (R\times \Sigma ),
d(e^a\lambda )\ominus d(e^b\lambda ))$$
and
$${\cal {L}}=\{ ((0,\sigma ),(0,\sigma ))|(0,\sigma )\in R\times \Sigma \}.$$
Let
$$L'={\cal {L}}\times R,
L_s'={\cal {L}}\times \{s\}.$$
Then
define
\begin{eqnarray}
&&G':L'\to V'\cr
&&G'(l')=G'(((\sigma ,0),(\sigma ,0)),s)
=((0,\sigma ),(0,\eta _s(\sigma )))
\end{eqnarray}
Then

$$W'=G'(L')=\{ ((0,\sigma ),(0,\eta _s(\sigma )))
|(0,\sigma )\in R\times \Sigma ,s\in R\}$$

$$W_s'=G'(L'_s)=\{ ((0,\sigma ),(0,\eta _s(\sigma )))
|(0,\sigma )\in R\times \Sigma \}$$
for fixed $s\in R$.

\begin{Lemma}
There does not exist any Reeb closed orbit in
$(\Sigma ,\lambda )$ if and only if
$W'_s\cap W'_{s'}$ is empty for $s\ne s'$.
\end{Lemma}
Proof. First if there exists a closed Reeb orbit in
$(\Sigma ,\lambda )$, i.e., there exists
$\sigma _0\in \Sigma $, $t_0>0$ such that
$\sigma _0=\eta _{t_0}(\sigma _0)$, then
$((0,\sigma _0),(0,\sigma _0))\in W'_0\cap W'_{t_0}$.
Second if there exists $s_0\ne s_0'$ such
that $W'_{s_0}\cap W'_{s_0'}\ne \emptyset $, i.e.,
there exists $\sigma _0$ such that
$$((0,\sigma _0),(0,\eta _{s_0}(\sigma _0))
=((0,\sigma _0),(0,\eta _{s_0'}(\sigma _0)),$$
then $\eta _{(s_0-s_0')}(\sigma _0)=\sigma _0$, i.e.,
$\eta _t(\sigma _0)$ is a closed Reeb orbit.

\begin{Lemma}
If there does not exist any closed Reeb orbit in
$(\Sigma ,\lambda )$ then
there exists a smooth Lagrangian injective immersion
$G':W'\to V'$ with $G'(((0,\sigma ),(0,\sigma )),s)
=((0,\sigma ),(0,\eta _s(\sigma )))$
such that
\begin{equation}
G'_{s_1,s_2}:{\cal L}\times (-s_1,s_2)\to V'
\end{equation}
is a regular exact Lagrangian embedding for any finite real number
$s_1$, $s_2$, here we denote by $W'(s_1,s_2)=G'_{s_1,s_2}({\cal L}\times
(s_1,s_2))$.
\end{Lemma}
Proof. One check
\begin{equation}
{G'}^*((e^a\lambda -e^b\lambda ))
=\lambda -\eta (\cdot ,\cdot )^*\lambda
=\lambda -(\eta _s^*\lambda +i_X\lambda ds)=-ds
\end{equation}
since $\eta _s^*\lambda =\lambda $.
This implies that ${G}'$ is an exact  Lagrangian embedding, this proves
Lemma 3.2.

\vskip 3pt

Now set
\begin{equation}
c(s,t)=\varepsilon te^{-s^2} \label{eq:fu}
\end{equation}
\begin{equation}
\psi _0(s,t)=se^{c(s,t)}c_s=-2e^{(\varepsilon te^{-s^2})-s^2}\varepsilon
ts^2 =\varepsilon \psi _0'\label{eq:fu0}
\end{equation}
here $\psi _0'(s,t)=-2ts^2e^{(\varepsilon te^{-s^2})-s^2}$;
\begin{equation}
\psi _1(s,t)=\int _{-\infty }^s\psi _0(\tau, t)d\tau
=\varepsilon \int _{-\infty }^s\psi _0'(\tau ,t)d\tau
=\varepsilon \psi '_1\label{eq:fu1}
\end{equation}
here $\psi _1'=\int _{-\infty }^s\psi _0'$;
\begin{equation}
\psi (s,t)={{\partial \psi _1}\over {\partial t}}-se^{c(s,t)}c_t
=\varepsilon \psi '\label{eq:fu2}
\end{equation}
here $\psi '(s,t)=
{{\partial \psi _1}\over {\partial t}}-se^{\varepsilon te^{-s^2}}
e^{-s^2}$;
\begin{equation}
\Psi '=se^{c(s,t)}; \tilde {l}'=-\psi (s,t)dt.  \label{eq:fu3}
\end{equation}

Now we construct an isotopy of Lagrangian injective immersions as follows:
\begin{eqnarray}
&&F':{\cal {L}}\times R\times [0,1]\to (R\times \Sigma )\times
(R\times \Sigma )\cr
&&F'(((0,\sigma ),(0,\sigma )),s,t)=((c(s,t),\sigma ),(c(s,t),\eta _s(\sigma ))) \cr
&&F'_t(((0,\sigma ),(0,\sigma )),s)=F'(((0,\sigma ),(0,\sigma )),s,t)
\end{eqnarray}
\begin{Lemma}
If there does not exist any closed Reeb orbit
in $(\Sigma ,\lambda )$ then for the choice of $c(s,t)$ satisfying
$\int _0^sc(s,t)ds$ or $\int _s^0c(s,t)ds $ exists and
is smooth on $(s,t)$, $F'$ is an exact isotopy of Lagrangian
embeddings. Moreover if $c(s,0)\ne c(s,1)$, then
$F'_0(\Sigma \times R)\cap F'_1(\Sigma \times R)=\emptyset $.
\end{Lemma}
Proof. Let $F'_t=F'(\cdot,t):{\cal {L}}\times R\to
(R\times \Sigma )\times (R\times \Sigma )$.
It is obvious that $F'_t$ is an embedding. We check
that
\begin{eqnarray}
{F'}^*(e^a\lambda \ominus e^b\lambda )&=&-e^{c(s,t)}ds \cr
&=&-\{ d(se^{c(s,t)} )-sde^{c(s,t)}\} \cr
&=&-\{ d(se^{c(s,t)})-se^{c(s,t)}c_sds-se^{c(s,t)}c_tdt \} \cr
&=&-\{ d(se^{c(s,t)})-d_s{\psi _1}-se^{c(s,t)}c_tdt \} \cr
&=&-\{ d((se^{c(s,t)})-\psi _1)
+{{\partial \psi _1}\over {\partial t}}dt
-se^{c(s,t)}c_tdt  \} \cr
&=&-\{d\Psi '+{{\partial \psi _1}\over {\partial t}}dt-se^{c(s,t)}c_tdt \}\cr
&=&-d\Psi '-\psi (s,t)dt \cr
&=&-d\Psi '+\tilde {l}' \label{eq:3.20}
\end{eqnarray}
here $\psi _i,\psi _i'$, and $\tilde {l}',\Psi'$ as in
(\ref{eq:fu0}-\ref{eq:fu3}).

   Let $(V',\omega ')$, $W'$ as above and
$(V,\omega )=(V'\times C,\omega '\oplus \omega _0)$.
As in \cite[p330,2.3.$B'_3$]{gro}(see also
\cite[p291-292]{aul}), we use figure eight trick invented by Gromov to
construct a Lagrangian submanifold in $V$ through the
Lagrange isotopy $F'$ in $V'$.
Fix a positive $\delta <1$ and take a $C^{\infty }$-map $\rho :S^1\to
[0,1]$, where the circle $S^1$ ia parametrized by $\Theta \in [-1,1]$,
such that the $\delta -$neighborhood $I_0$ of $0\in S^1$ goes to
$0\in [0,1]$ and $\delta -$neighbourhood $I_1$ of $\pm 1\in S^1$
goes $1\in [0,1]$. Let
\begin{eqnarray}
\tilde {l}&=&\psi (s,\rho (\Theta ))\rho '(\Theta )d\Theta \cr
&=&\Phi d\Theta
\end{eqnarray}
be the pull-back of the form
$\tilde {l}'=\psi (s,t)dt $ to $W'\times S^1$ under the map
$(w',\Theta )\to (w',\rho (\Theta ))$ and
assume without loss of generality $\Phi $
vanishes on $W'\times (I_0\cup I_1)$.

  Next, consider a map $\alpha $ of the annulus $S^1\times [\Phi _-,\Phi _+]$
into $R^2$, where $\Phi _-$ and $\Phi _+$ are the lower and the upper
bound of the fuction $\Phi $ correspondingly, such that

   $(i)$ The pull-back under $\alpha $ of the form
$dx\wedge dy$ on $R^2$ equals $-d\Phi \wedge d\Theta $.

   $(ii)$ The map $\alpha $ is bijective on $I\times [\Phi _-,\Phi _+]$
where $I\subset S^1$ is some closed subset,
such that $I\cup I_0\cup I_1=S^1$; furthermore, the origin
$0\in R^2$ is a unique double point of the map $\alpha $ on
$S^1\times 0$, that is
$$0=\alpha (0,0)=\alpha (\pm 1,0),$$
and
$\alpha $ is injective on $S^1=S^1\times 0$ minus $\{ 0,\pm 1\}$.

   $(iii)$ The curve $S^1_0=\alpha (S^1\times 0)\subset R^2$ ``bounds''
zero area in $R^2$, that is $\int _{S^1_0}xdy=0$, for the $1-$form
$xdy$ on $R^2$.
\begin{Proposition}
Let $V'$, $W'$ and $F'$ as above. Then there exists
an exact Lagrangian embedding $F:W'\times S^1\to V'\times R^2$
given by $F(w',\Theta )=(F'(w',\rho (\Theta )),\alpha (\Theta ,\Phi ))$.
\end{Proposition}
Proof. We follow as in \cite[2.3$B_3'$]{gro}. Now let
$F^*:W'\times S^1\to V'\times R^2$ be given by
$(w',\Theta )\to (F'(w,\rho (\Theta )),\alpha (\Theta ,\Phi ))$.
Then

   $(i)'$ The pull-back under $F^*$ of the form
$\omega =\omega '+dx\wedge dy$ equals
$d\tilde {l}^*-d\Phi \wedge d\Theta =0$ on $W'\times S^1$.

   $(ii)'$ The set of double points of $F^*$ is
$W'_0\cap W_1'\subset V'=V'\times 0\subset V'\times R^2$.

   $(iii)'$  If $F^*$ has no double point then the
Lagrangian submanifold $W=F^*(W'\times S^1)\subset
(V'\times R^2,\omega '+dx\wedge dy)$ is exact if and only if
$W_0'\subset V'$ is such.

   This completes the proof of Proposition 3.1.

\subsection{Formulation of Hilbert bundles}

Let $(\Sigma ,\lambda )$ be a closed $(2n-1)-$ dimensional manifold
with a contact form $\lambda $.
Let $S\Sigma =R\times \Sigma $ and put
$\xi =\ker (\lambda )$.
Let $J'_\lambda $ be
an almost complex structure
on $S\Sigma $
tamed by
the symplectic form $d(e^a\lambda )$.

We define a metric $g_\lambda $ on
$S\Sigma =R\times \Sigma $ by
\begin{equation}
g_\lambda =d(e^a\lambda )(\cdot ,J_\lambda \cdot )
\end{equation}
which is adapted to $J_\lambda $ and $d(e^a\lambda )$ but not
complete.

\vskip 3pt

In the following we denote by $(V',\omega ')=
((R\times \Sigma )\times (R\times \Sigma ),d(e^a\lambda _1-e^b\lambda _2))$
and
$(V,\omega )=(V'\times R^2,\omega '+dx\wedge dy)$
with the metric $g=g'\oplus g_0=g_{\lambda _1}
\oplus g_{\lambda _2}\oplus g_0 $ induced by
$\omega (\cdot ,J\cdot )$($J=J'\oplus i=J_{\lambda _1}
\oplus (-J_{\lambda _2})\oplus i$ and
$W\subset V$ a Lagrangian submanifold which was constructed in
section 3.1.

   Let $\bar V=D\times V$, then $\pi _1:\bar V\to D$ be
a symplectic vector bundle. Let $\bar J$ be an almost
complex structure on $\bar V$ such that $\pi _1:\bar V\to D$ is
a holomorphic map and each fibre $\bar V_z=\pi _1(z)$ is
a $\bar J$ complex submanifold. Let
$H^k(D)$ be the space of $H^k-$maps
from $D$ to $\bar V$, here
$H^k$ represents Sobolev derivatives up to order $k$. Let
$\bar W=\partial D\times W$, $\bar p=\{1\}\times p$,
$W^{\pm}=\{\pm i\}\times W $ and

$${\cal D}^k
=\{ \bar u \in H^k(D)|
\bar u(x)\in \bar W \ a.e \ for \ x\in \partial D \ and \ \bar u(1)=\bar p,
\bar u(\pm i)\in \{\pm i\}\times W\}$$
for $k\geq 100$.
\begin{Lemma}
Let $W$ be a closed Lagrangian submanifold in
$V$. Then,
${\cal D}^k$
is a pseudo-Hilbert manifold with the tangent bundle
\begin{equation}
T{\cal D}^k
=\bigcup _{\bar u\in {\cal {D}}^k}
\Lambda ^{k-1}
\end{equation}
here
$$\Lambda ^{k-1}=\{ \bar w\in H^{k-1}(\bar u^*(T\bar V)|
\bar w(1)=0, and \ \bar w(\pm i)\in T\bar W\} $$
\end{Lemma}

\begin{Note}
Since $W$ is not regular we know that ${\cal D}^k$ is
in general complete, however it is enough for our purpose.
\end{Note}
Proof: See \cite[p309-310]{aul} or follow step by step
from \cite[ch1]{kl}.

\vskip 3pt

   Now we consider  a section
from ${\cal D}^k$ to
$T{\cal D}^k$ follows as in
\cite[p327,2.2]{gro} or \cite[p310]{gro}, i.e.,
let $\bar \partial :{\cal D}^k\rightarrow T{\cal D}^k$
be the Cauchy-Riemmann section
\begin{equation}
\bar \partial \bar u={{\partial \bar u}\over {\partial s}}
+J{{\partial \bar u}\over {\partial t}}  \label{eq:CR}
\end{equation}
for $\bar u\in {\cal D}^k$.

\begin{Theorem}
The Cauchy-Riemann section $\bar \partial $ defined in (\ref{eq:CR})
is a Fredholm section of Index zero.
\end{Theorem}
Proof. According to the definition of the Fredholm section,
we need to prove that
$\bar u\in {\cal D}^k$, the linearization
$D\bar \partial (\bar u)$ of $\bar \partial $ at $\bar u$ is
a linear Fredholm
operator.
Note that
\begin{equation}
D\bar \partial (\bar u)=D{\bar \partial _{[\bar u]}}
\end{equation}
where
\begin{equation}
(D\bar \partial _{[\bar u]})v=\frac{\partial \bar v}{\partial s}
+J\frac{\partial \bar v}{\partial t}+A(\bar u)\bar v
\end{equation}
with
$$\bar v|_{\partial D}\in (\bar u|_{\partial D})^*T\bar W$$
here $A(\bar u)$ is $2n\times 2n$
matrix induced by the torsion of
almost complex structure, see \cite[p324,2.1]{gro} for the computation.

   Observe that the linearization $D\bar \partial (\bar u)$ of
$\bar \partial $ at $\bar u$ is equivalent to the following Lagrangian
boundary value problem
\begin{eqnarray}
&&{{\partial \bar v}\over {\partial s}}+\bar J
{{\partial \bar v}\over {\partial t}}
+A(\bar u)\bar v=\bar f, \ \bar v\in \Lambda ^k(\bar u^*T\bar V)\cr
&&\bar v(t)\in T_{\bar u(t)}W, \ \ t\in {\partial D}  \label{eq:Lin}
\end{eqnarray}
One
can check that (\ref{eq:Lin})
defines a linear Fredholm operator. In fact,
by proposition 2.2 and Lemma 2.1, since the operator $A(\bar u)$ is a compact,
we know that the operator $\bar \partial $ is a nonlinear Fredholm operator
of the index zero.

\begin{Definition}
Let $X$ be a Banach manifold and $P:Y\to X$ the Banach
vector bundle.
A Fredholm section $F:X\rightarrow Y$ is
proper if $F^{-1}(0)$ is a compact set and is called
generic if $F$ intersects the zero section transversally, see \cite[p327-328,2.2B]{gro}.
\end{Definition}
\begin{Definition}
$deg(F,y)=\sharp \{ F^{-1}(0)\} mod2$ is called the Fredholm
degree of a Fredholm section (see\cite[p327-328,2.2B]{gro}).
\end{Definition}
\begin{Theorem}
Assum that $\bar J=i\oplus J$ on $\bar V$ and $i$ is complex structure
on $D$ and $J$ the almost complex structure on
$V$. Assume that $J$ is integrable at $p\in V$. Then the Fredholm section
$F=\bar \partial _{\bar J}: {\cal D}^k\rightarrow T{\cal D}^k$
constructed in (\ref{eq:CR}) has degree one, i.e.,
$$deg(F,0)=1$$
\end{Theorem}
Proof: We assume that $\bar u:D\mapsto \bar V$ be a $\bar J-$holomorphic disk
with boundary $\bar u(\partial D)\subset \bar W$ and
by the assumption that $\bar u$ is homotopic to the
map $\bar u_1=(id,p)$. Since almost complex
structure ${\bar J}$ splits and
is tamed by  the symplectic form $\bar \omega $,
by stokes formula,
we conclude the second component $u: D\rightarrow
V$ is a constant map. Because $u(1)=p$, We know that
$F^{-1}(0)=(id,p)$.
Next we show that the linearizatioon $DF_{(id,p)}$ of $F$ at $(id,p)$ is
an isomorphism from $T_{(id,p)}{\cal D}^k$ to $E$.
This is equivalent to solve the equations
\begin{eqnarray}
{\frac {\partial \bar v}{\partial s}}+\bar J{\frac {\partial \bar v}{\partial t}}
=f\\
\bar v|_{\partial D}\subset T_{(id ,p)}\bar W
\end{eqnarray}
here $\bar J=i+J(p)$ since $J$ is integrable at $p$. By Lemma 2.1,
we know that $DF_{(id,p)}$ is an isomorphism.
Therefore $deg(F,0)=1$.

\section{Anti-holomorphic sections}

In this section we construct a Fredholm
section which is not proper as in
\cite[p329-330,2.3.B]{gro}(see also \cite[p315, 5.3]{aul}).

  Let $(V',\omega ')=(S\Sigma \times S\Sigma ,
d(e^a\lambda _1-e^b\lambda _2))$ and
$(V,\omega )=(V'\times C, \omega '\oplus \omega _0)$,
$W$ as in section3 and $J=J'\oplus i$, $g=g'\oplus g_0$,
$g_0$ the standard metric on $C$.

   Now let $c\in C$ be a non-zero vector. We consider
$c$ as an anti-holomorphic homomorphism
$c:TD\to TV'\oplus TC$, i.e., $c({{\partial }\over
{\partial \bar z}})
=(0,c\cdot {{\partial }\over {\partial z}}).$
Since the constant section $c$ is not a section of the
Hilbert bundle in section 3 due to $c$ is not
tangent to the Lagrangian submanifold $W$, we must modify it as follows:

\vskip 3pt

  Let $c$ as in section 4.1, we define
\begin{eqnarray}
c_{\chi ,\delta }(z,v)=\left\{ \begin{array}{ll}
c \ \ \ &\mbox{if\  $|z|\leq 1-2\delta $,}\cr
0 \ \ \ &\mbox{otherwise}
\end{array}
\right.
\end{eqnarray}
Then by using the cut off function $\varphi _h(z)$ and
its convolution with section
$c_{\chi ,\delta }$, we obtain a smooth section
$c_\delta$ satisfying

\begin{eqnarray}
&&c_{\delta }(z,v)=\left\{ \begin{array}{ll}
c \ \ \ &\mbox{if\  $|z|\leq 1-3\delta $,}\cr
0 \ \ \ &\mbox{if\  $|z|\geq 1-\delta $.}
\end{array}
\right.   \cr
&&|c_\delta |\leq |c|
\end{eqnarray}
for $h$ small enough, for the convolution theory see
\cite[ch1,p16-17,Th1.3.1]{hor}.
Then one can easily check
that $\bar c_\delta =(0,0,c_\delta )$
is an anti-holomorphic section tangent to $\bar W$.

\vskip 3pt

Now we modify the almost conplex structure
on the $V$. Let $J_1$, $J_2$ be the almost
complex structures on $V$ tamed by $\omega $. Let
$g_i=\omega (\cdot ,J_i\cdot )$ the metrics
by $\omega $ and $J_i$. We assume there
exists a constant $c_1$ such that
\begin{equation}
c_1^{-1}g_1\leq g_2\leq c_1g_1
\end{equation}
Let
\begin{eqnarray}
J_{\chi ,\delta }(z,v)=\left\{ \begin{array}{ll}
i\oplus J_1\ \ \ &\mbox{if\  $|z|\leq 1-2\delta $,}\cr
i\oplus J_2 \ \ \ &\mbox{otherwise}
\end{array}
\right.
\end{eqnarray}
Then by using the cut off function $\varphi _h(z)$ and
its convolution with section
$J_{\chi ,\delta }$, we obtain a smooth section
$J_\delta$ satisfying

\begin{eqnarray}
J_{\delta }(z,v)=\left\{ \begin{array}{ll}
i\oplus J_1\ \ \ &\mbox{if\  $|z|\leq 1-3\delta $,}\cr
i\oplus J_2\ \ \ &\mbox{if\  $|z|\geq 1-\delta $.}
\end{array}
\right.
\end{eqnarray}
for $h$ small enough, for the convolution theory see
\cite[ch1,p16-17,Th1.3.1]{hor}.

Now we get an almost conplex structure
$\bar J=i\oplus J_\delta $
on the symplectic fibration $D\times V \to D$
such that $\pi _1:D\times V\to D$ is a holomorphic
fibration and $\pi _1^{-1}(z)$ is an almost complex submanifold.
Let $g_\delta =\bar \omega (\cdot ,\bar J\cdot ) $,
$\bar g_i=g_0\oplus g_i$ be the metrics on
$D\times V$, $g_0$ is metric on $D$.
We assume there exists a constants $c_2$ such that
\begin{equation}
c_2^{-1}\bar g_i\leq g_\delta \leq c_2g_i, i=1,2.
\end{equation}

Now we consider the
equations
\begin{eqnarray}
&&\bar v=(id ,v)=(id, v',f):D\to D\times V'\times C \cr
&&\bar \partial _{J_\delta }v=c_\delta \cr
&&\bar \partial _{J'}v'=0,\bar \partial f=c_\delta
\ on \ D_{1-2\delta }\cr
&&v|_{\partial D}:\partial D\to W\ \ \  \label{eq:4.16}
\end{eqnarray}
here $v$ homotopic to constant map
$\{ p\}$ relative to $W$.
Note that $W\subset V\times B_{R}(0)$ for $2\pi R(\varepsilon )^2$,
here $R(\varepsilon )\to 0$ as $\varepsilon \to 0$ and
$\varepsilon $ as in section 3.1.

\begin{Lemma}
Let $\bar v=(id ,v)$ be the solutions of (\ref{eq:4.16}), then one has
the following estimates
\begin{eqnarray}
E({v})=
\{
\int _D(g'({{\partial {v'}}\over {\partial x}},
{J'}{{\partial {v'}}\over {\partial x}})
+g'({{\partial {v'}}\over {\partial y}},
{J'}{{\partial {v'}}\over {\partial y}}) \nonumber \\
+g_0({{\partial {f}}\over {\partial x}},
{i}{{\partial {f}}\over {\partial x}})
+g_0({{\partial {f}}\over {\partial y}},
{i}{{\partial {f}}\over {\partial y}}))d\sigma \}
\leq 4\pi R(\varepsilon )^2.
\end{eqnarray}
\end{Lemma}
Proof: Since $v(z)=(v'(z),f(z))$ satisfy (\ref{eq:4.16})
and $v(z)=(v'(z),f(z))\in V'\times C$
is homotopic to constant map $v_0:D\to \{ p\}\subset W$
in $(V,W)$, by the Stokes formula
\begin{equation}
\int _{D}v^*(\omega '\oplus \omega _0)=0
\end{equation}
Note that the metric $g$ is adapted to the symplectic form
$\omega $ and $J$, i.e.,
\begin{equation}
g=\omega  (\cdot ,J\cdot )
\end{equation}
By the simple algebraic computation, we have
\begin{equation}
\int _{D}{v}^*\omega  ={{1}\over {4}}
\int _{D^2}(|\partial v|^2
-|\bar {\partial }v|^2)=0
\end{equation}
and
\begin{equation}
|\nabla v|={{1}\over {2}}(
|\partial v|^2 +|\bar \partial v|^2
\end{equation}
Then
\begin{eqnarray}
E(v)&=&\int _{D} |\nabla v| \nonumber \\
      &=&\int _{D}\{ {{1}\over {2}}(
|\partial v|^2+|\bar \partial v|^2)\})d\sigma \nonumber \\
&=&\int _D|c_\delta |_{\bar g}^2d\sigma
\end{eqnarray}
By Cauchy integral formula,
\begin{equation}
f(z)={{1}\over {2\pi i}}\int _{\partial D}{{f(\xi)}\over
{\xi -z}}d\xi
+{{1}\over {2\pi i}}
\int _D{{\bar \partial f(\xi )}\over {\xi -z}}d\xi\wedge d\bar \xi
\label{eq:Ch}
\end{equation}
Since $f$ is smooth up to the boundary,
we integrate the two sides on $D_r$ for $r<1$, one get
\begin{eqnarray}
\int _{\partial D_r}f(z)dz&=&
\int _{\partial D_r}{{1}\over {2\pi i}}\int _{\partial D}{{f(\xi)}\over
{\xi -z}}d\xi dz
+\int _{\partial D_r}{{1}\over {2\pi i}}
\int _D{{\bar \partial f(\xi )}\over {\xi -z}}d\xi\wedge d\bar \xi \cr
&=&0+{{1}\over {2\pi i}}
\int _D\int _{\partial D_r}{{\bar \partial
f(\xi )}\over {\xi -z}}dzd\xi\wedge d\bar \xi \cr
&=&{{1}\over {2\pi i}}
\int _D2\pi i\bar \partial
f(\xi )d\xi\wedge d\bar \xi \label{eq:Ch1}
\end{eqnarray}
Let $r\to 1$, we get

\begin{eqnarray}
\int _{\partial D}f(z)dz=
\int _D\bar \partial
f(\xi )d\xi\wedge d\bar \xi \label{eq:Ch2}
\end{eqnarray}

By the equations (\ref{eq:4.16}),
one get
\begin{equation}
\bar \partial f=c \ on \ D_{1-2\delta }
\end{equation}
So, we have
\begin{equation}
2\pi i(1-2\delta )c=\int _{\partial D}f(z)dz-
\int _{D-D_{1-2\delta }}\bar \partial
f(\xi )d\xi\wedge d\bar \xi \label{eq:Ch3}
\end{equation}
So,
\begin{eqnarray}
|c|&\leq &{{1}\over {2\pi (1-2\delta )}}|\int _{\partial D}f(z)dz|
+|\int _{D-D_{1-2\delta }}\bar \partial
f(\xi )d\xi\wedge d\bar \xi |\cr
&\leq &{{1}\over {2\pi (1-2\delta )}}
2\pi |diam(pr_2(W))+c_1c_2|c|(\pi -\pi (1-2\delta )^2))
\label{eq:Ch4}
\end{eqnarray}
Therefore, one has
\begin{eqnarray}
|c|&\leq &
c(\delta )R(\varepsilon )\label{eq:Ch5}
\end{eqnarray}
and
\begin{eqnarray}
E(v)&=&\pi \int _D|c_\delta |_{\bar g}^2  \cr
    &=&\pi c(\delta )^2R(\varepsilon )^2.
\end{eqnarray}
This finishes the proof of Lemma.

\begin{Proposition}
For $|c|\geq 2c(\delta )R(\varepsilon )$, then the
equations (\ref{eq:4.16})
has no solutions.
\end{Proposition}
Proof. By \ref{eq:Ch5}, it is obvious.

\begin{Theorem}
The Fredholm section $F_1=\bar \partial _{\bar J}+\bar c_\delta
: {\cal  {D}}^k\rightarrow E$ is not proper.
\end{Theorem}
Proof. By the Proposition 4.1 and Theorem 3.2,
it is obvious(see also \cite[p330,$2.3B_1$]{gro} or
\cite[p316]{aul}).

\section{$J-$holomorphic section}

Recall that $W(-K,K)\subset W\subset V'\times R^{2}$ as in section 3.
The Riemann metric $g$ on $V'\times R^{2}$
induces a metric $g|W$.

   Now let $c\in C$ be a non-zero vector and
$c_\delta $ the induced anti-holomorphic section. We consider the
nonlinear inhomogeneous equations (\ref{eq:4.16}) and
transform it into $\bar J-$holomorphic map by
considering its graph as in \cite[p319,1.4.C]{gro} or
\cite[p312,Lemma5.2.3]{aul}.

Denote by $Y^{(1)}\to D\times V$ the bundle of homomorphisms $T_s(D)\to
T_v(V)$. If $D$ and $V$ are given the disk and the almost
K\"ahler manifold, then
we distinguish the subbundle $X^{(1)}\subset Y^{(1)}$ which consists of
complex linear homomorphisms and we denote $\bar X^{(1)}\to D\times V$ the
quotient bundle $Y^{(1)}/X^{(1)}$. Now, we assign to each $C^1$-map $
v:D\to V$ the section $\bar \partial v$ of the bundle $\bar X^{(1)}$ over
the graph $\Gamma _v\subset D\times V$ by composing the differential of $v$
with the quotient homomorphism $Y^{(1)}\to \bar {X}^{(1)}$. If $c_\delta
:D\times
V\to \bar X$ is a $H^k-$ section we write $\bar \partial v=c_\delta $
for the
equation $\bar \partial v=c_\delta |\Gamma _v$.

\begin{Lemma}
(Gromov\cite[1.4.$C'$]{gro})There exists a unique almost complex
structure $J_g$ on $D\times V$(which
also depends on the given structures in $D$ and in $V$), such that
the (germs of) $J_\delta-$holomorphic sections $v:D\to D\times V$ are exactly and
only the solutions
of the equations $\bar \partial v=c_\delta $. Furthermore, the
fibres $z\times V\subset D\times V$ are $J_\delta-$holomorphic(
i.e. the subbundles $T(z\times V)\subset T(D\times V)$ are $J_\delta-$complex)
and the structure
$J_\delta|z\times V$ equals the original structure on $V=z\times V$.
Moreover $J_\delta $ is tamed by $k\omega _0\oplus \omega $ for
$k$ large enough which is independent of $\delta $.
\end{Lemma}

\section{Gromov's $C^0-$convergence theorem}

\subsection{Analysis of Gromov's figure eight}

Since $W'\subset S\Sigma \times S\Sigma $ is an exact Lagrangian
submanifold and $F'_t$ is an exact Lagrangian isotopy(see section
3.1). Now we carefully check the Gromov's construction of
Lagrangian submanifold $W\subset V'\times R^2$ from the exact
Lagrangian isotopy of $W'$ in section 3.

   Let $S^1\subset T^*S^1$ be a zero section and $S^1=\cup _{i=1}^4S_i$
be a partition of the zero section $S^1$ such that $S_1=I_0$,
$S_3=I_1$. Write $S^1\setminus \{I_0\cup I_1\}=I_2\cup I_3$ and
$I_0=(-\delta ,-{{5}\over {6}}\delta ]\cup (-{{5\delta }\over
{6}},+{{5\delta }\over {6}})\cup [{{5\delta }\over {6}},
\delta)=I_0^-\cup I_0'\cup I_0^+$, similarly $I_1=(1-\delta
,1-{{5}\over {6}}\delta ]\cup (1-{{5\delta }\over {6}},1+{{5\delta
}\over {6}})\cup [1+{{5\delta }\over {6}}, 1+\delta)=I_1^-\cup
I_1'\cup I_1^+$. Let $S_2=I_0^+\cup I_2\cup I_1^-$, $S_4=I_1^+\cup
I_3\cup I_0^+$. Moreover, we can assume that the double points of
map $\alpha $ in Gromov's figure eight is contained in $(\bar
I_0'\cup \bar I_1')\times [\Phi _-,\Phi _+]$, here $\bar
I_0'=(-{{5\delta }\over {12}},+{{5\delta }\over {12}})$ and $\bar
I_1'=(1-{{5\delta }\over {12}},1+{{5\delta }\over {12}})$. Recall
that $\alpha : (S^1\times [5\Phi _-,5\Phi _+])\to R^2$ is an exact
symplectic immersion, i.e., $\alpha ^*(-ydx)-\Psi d\Theta =dh$,
$h:T^*S^1\to R$. By the construction of figure eight, we can
assume that $\alpha '_i=\alpha |((S^1\setminus I'_i)\times [5\Phi
_-,5\Phi _+])$ is an embedding for $i=0,1$.
 Let $Y=\alpha (S^1\times
[5\Phi _-,5\Phi _+])\subset R^2$ and $Y_i=\alpha (S_i\times [5\Phi
_-,5\Phi _+])\subset R^2$. Let $\alpha _i=\alpha |Y_i(S^1\times
[5\Phi _-,5\Phi _+])$. So, $\alpha _{i}$ puts the function $h$ to
the function $h_{i0}=\alpha _{i}^{-1*}h$ on $Y_i$. We extend the
function $h_{i0}$ to whole plane $R^2$. In the following we take
the liouville form $\beta _{i0}=-ydx -dh_{i0} $ on $R^2$. This
does not change the symplectic form $dx\wedge dy $ on $R^2$. But
we have $\alpha _i^*\beta =\Phi d\Theta $ on $(S_i\times [5\Phi
_-,5\Phi _+])$ for $i=1,2,3,4$.  Finally, note that
\begin{eqnarray}
&&F:W'\times S^1\to V'\times R^2;\cr &&F(w',\Theta )=(F'_{\rho
(\Theta )}(w'),\alpha (\Theta ,\Phi (w',\rho (\Theta )).
\end{eqnarray}
Since $\rho (\Theta )=0 $ for $\Theta \in I_0$ and $\rho (\Theta
)=1$ for $\Theta \in I_1$, we know that $\Phi (w',\rho (\Theta
))=0$ for $\Theta \in I_0\cup I_1$. Therefore,
\begin{eqnarray}
F(W'\times I_0)=W'\times \alpha (I_0); F(W'\times I_1)=W'\times
\alpha (I_1).
\end{eqnarray}

\subsection{Gromov's Schwartz lemma}

In our proof we need a crucial tools, i.e., Gromov's Schwartz
Lemma as in \cite{gro}. We first consider the case without
boundary.
\begin{Proposition}
Let $(V,J,\mu )$ be as in section 4 and $V_K$ the compact part of
$V$.  There exist constants $\varepsilon _0$ and $C$(depending
only on the $C^0-$ norm of $\mu $ and on the $C^\alpha $ norm of
$J$ and $A_0$) such that every $J-$holomorphic map of the unit
disc to an $\varepsilon _0$-ball of $V$ with center in $V_K$ and
area less than $A_0$ has its derivatives up to order $k+1+\alpha $
on $D_{{1}\over {2}}(0)$ bounded by $C$.
\end{Proposition}
For a proof, see\cite{gro}.

\vskip 3pt

Now we consider the Gromov's Schwartz Lemma for $J-$holomorphic
map with boundary in a closed Lagrangian submanifold as in
\cite{gro}.

\begin{Proposition}
Let $(V,J,\mu )$ as above and $L\subset V$ be a closed Lagrangian
submanifold and $V_K$ one compact part of $V$. There exist
constants $\varepsilon _0$ and $C$(depending only on the $C^0-$
norm of $\mu $ and on the $C^\alpha $ norm of $J$ and $K,A_0$)
such that every $J-$holomorphic map of the half unit disc $D^+$ to
a $\varepsilon _0$-ball of $V$ with boundary in $L$ and area less
than $A_0$ has its derivatives up to order $k+1+\alpha $ on
$D_{{1}\over {2}}^+(0)$ bounded by $C$.
\end{Proposition}
For a proof see \cite{gro}.

\vskip 3pt

Since in our case $W$ is a non-compact Lagrangian submanifold,
Proposition 6.2 can not be used directly but the proofs of
Proposition 6.1-2 is still holds in our case.

\begin{Lemma}
Recall that $V=V'\times R^2$. Let $(V,J,\mu )$ as above and
$W\subset V$ be as above and $V_c$ the compact set in $V$. Let
$\bar V=D\times V$, $\bar W=\partial D\times W$, and $\bar
V_c=D\times V_c$. Let $Y=\alpha (S^1\times [5\Phi _-,5\Phi
_+])\subset R^2$. Let $Y_i=\alpha (S_i\times [5\Phi _-,5\Phi
_+])\subset R^2$. Let $\{ X_j\}_{j=1}^q$ be a Darboux covering of
$\Sigma $ and $V'_{ij}=(R\times X_i)\times (R\times X_j)$.
 Let $\partial D=S^{1+}\cup S^{1-}$. There exist
constant $c_0$ such that every $J-$holomorphic map $v$ of the half
unit disc $D^+$ to the $D\times V_j'\times R^2$ with its boundary
$v((-1,1))\subset (S^{1\pm})\times F({\cal L}\times R\times
S_i)\subset \bar W,i=1,..,4$ has
\begin{equation}
area(v(D^+))\leq c_0l^2(v(\partial 'D^+)).
\end{equation}
here $\partial 'D^+=\partial D\setminus [-1,1]$ and $l(v(\partial
'D^+))=length(v(\partial 'D^+))$.
\end{Lemma}
Proof. Let $\bar W_{i\pm }=S^{1\pm}\times F(W'\times S_i)$. Let
$v=(v_1,v_2):D^+\to \bar V=D\times V$ be the $J-$holomorphic map
with $v(\partial D^+)\subset \bar W_{i\pm}\subset \partial D\times
W$, then

\begin{eqnarray}
area(v)&=&\int _{D^+}v^*d(\alpha _0\oplus \alpha )\cr
&=&\int_{D^+}dv^*(\alpha _0\oplus \alpha )\cr &=&\int _{\partial
D^+}v^*(\alpha _0 \oplus \alpha )\cr &=&\int _{\partial
D^+}v_1^*\alpha _0+\int _{\partial D^+}v_2^*\alpha \cr &=&\int
_{\partial 'D^+\cup [-1,+1]}v_1^*\alpha _0 +\int _{\partial
'D^+\cup [-1,+1]}v_2^*(e^a\lambda-ydx-dh_{i0})\cr &=&\int
_{\partial 'D^+\cup [-1,+1]}v_1^*\alpha _0 + \int _{\partial
'D^+}v_2^*(e^a\lambda-ydx-dh_{i0})+B_1 ,\label{eq:al1}
\end{eqnarray}
here $B_1=\int _{[-1,+1]}{v_2}^*(-d\Psi ')$. Now take a zig-zag
curve $C$ in $V_j'\times Y_i$ connecting $v_2(-1)$ and $v_2(+1)$
such that

\begin{eqnarray}
\int _{C} (e^a\lambda+ydx)&=&B_1\cr length(C)&\leq&
k_1length(v_2(\partial 'D^+)) \label{eq:m_1}
\end{eqnarray}
Now take a minimal surface $M$ in $V'_{ij}\times R^2$ bounded by
$v_2(\partial 'D^+)\cup C$, then by the isoperimetric
ineqality(see[\cite[p283]{grob}), we get

\begin{eqnarray}
area(M)&\leq &m_1length (C+v_2(\partial 'D^+))^2\cr &\leq &
m_2length (v_2(\partial 'D^+))^2,
\end{eqnarray}
here we use the (\ref{eq:m_1}).

Since $area(M)\geq \int _M\omega $ and $\int _M\omega =\int
_{D^+}v_2^*\omega =area(v)$, this proves the lemma.

\begin{Lemma}
Let $v$ as in Lemma 6.1, then we have
\begin{equation}
area(v(D^+)\geq c_0(dist(v(0),v(\partial 'D^+)))^2,
\end{equation}
here $c_0$ depends only on $\Sigma ,J,\omega,...,$etc, not on $v$.
\end{Lemma}
Proof. By the standard argument as in \cite[p79]{aul}.

\vskip 3pt

The following estimates is a crucial step in our proof.

\begin{Lemma}
Recall that $V=V'\times R^2$. Let $(V,J,\mu )$ as above and
$W\subset V$ be as above and $V_c$ the compact set in $V$. Let
$\bar V=D\times V$, $\bar W=\partial D\times W$, and $\bar
V_c=D\times V_c$. Let $Y=\alpha (S^1\times [5\Phi _-,5\Phi
_+])\subset R^2$. Let $Y_i=\alpha (S_i\times [5\Phi _-,5\Phi
_+])\subset R^2$. Let $\partial D=S^{1+}\cup S^{1-}$. There exist
constant $c_0$ such that every $J-$holomorphic map $v$ of the half
unit disc $D^+$ to the $D\times V'\times R^2$ with its boundary
$v((-1,1))\subset (S^{1\pm})\times F({\cal L}\times R\times
S_i)\subset \bar W,i=1,..,4$ has
\begin{equation}
area(v(D^+))\leq c_0l^2(v(\partial 'D^+)).
\end{equation}
here $\partial 'D^+=\partial D\setminus [-1,1]$ and $l(v(\partial
'D^+))=length(v(\partial 'D^+))$.
\end{Lemma}
Proof. We first assume that $\varepsilon $ in section 3.1 is small
enough. Let $l_0$ is a constant small enough. If $length(\partial
'D^+)\geq l_0$, then Lemma 6.3 holds. If $length(\partial
'D^+)\leq l_0$ and $v(D^+)\subset D\times V'_{ij}\times R^2$, then
Lemma6.3 reduces to Lemma6.1. If $length(\partial 'D^+)\leq l_0$
and $v(D^+)\bar \subset D\times V'_{ij}\times R^2$, then Lemma6.2
imples $area(v)\geq \tau _0>100\pi R(\varepsilon )^2$, this is a
contradiction. Therefore we  proved the lemma.

\begin{Proposition}
Let $(V,J,\mu )$ and $W\subset V$ be as in section 4 and $V_K$ the
compact part of $V$. Let $\bar V$, $\bar V_K$ and $\bar W$ as
section 5.1. There exist constants $\varepsilon _0$ (depending
only on the $C^0-$ norm of $\mu $ and on the $C^\alpha $ norm of
$J$) and $C$(depending only on the $C^0$ norm of $\mu $ and on the
$C^{k+\alpha }$ norm of $J$) such that every $J-$holomorphic map
of the half unit disc $D^+$ to the $D\times V'\times R^2$ with its
boundary $v((-1,1))\subset (S^{1\pm})\times F({\cal L}\times
R\times S_i)\subset \bar W,i=1,..,4$ has its derivatives up to
order $k+1+\alpha $ on $D_{{1}\over {2}}^+(0)$ bounded by $C$.
\end{Proposition}
Proof. One uses Lemma 6.3 and Gromov's proof on Schwartz lemma to
yield proposition 6.3.

\subsection{Removal singularity of $J-$curves}

In our proof we need another crucial tools, i.e., Gromov's removal
singularity theorem\cite{gro}. We first consider the case without
boundary.
\begin{Proposition}
Let $(V,J,\mu )$ be as in section 4 and $V_K$ the compact part of
$V$. If $v:D\setminus \{0\}\to V_K$ be a $J-$holomorphic disk with
bounded energy and bounded image, then $v$ extends to a
$J-$holomorphic map from the unit disc $D$ to $V_K$.
\end{Proposition}
For a proof, see\cite{gro}.

\vskip 3pt

Now we consider the Gromov's removal singularity theorem for
$J-$holomorphic map with boundary in a closed Lagrangian
submanifold as in \cite{gro}.

\begin{Proposition}
Let $(V,J,\mu )$ as above and $L\subset V$ be a closed Lagrangian
submanifold and $V_K$ one compact part of $V$. If $v:(D^+\setminus
\{0\},\partial ''D^+\setminus \{0\})\to (V_K,L)$ be a
$J-$holomorphic half-disk with bounded energy and bounded image,
then $v$ extends to a $J-$holomorphic map from the half unit disc
$(D^+,\partial ''D^+)$ to $(V_K,L)$.
\end{Proposition}
For a proof see \cite{gro}.

\vskip 3pt

\begin{Proposition}
Let $(V,J,\mu )$ and $W\subset V$ be as in section 4 and $V_c$ the
compact set in $V$. Let $\bar V=D\times V$, $\bar W=\partial
D\times W$, and $\bar V_c=D\times V_c$. Then every $J-$holomorphic
map $v$ of the half unit disc $D^+\setminus \{0\}$ to the $\bar V$
with center in $\bar V_c$ and its boundary $v((-1,1)\setminus
\{0\})\subset (S^{1\pm})\times F({\cal L}\times [-K,K]\times
S_i)\subset \bar W$ and
\begin{equation}
area(v(D^+\setminus \{0\}))\leq E
\end{equation}
extends to a $J-$holomorphic map $\tilde v:(D^+, \partial ''D)\to
(\bar V_c, \bar W)$.
\end{Proposition}
Proof. This is ordinary Gromov's removal singularity theorem by
$K-$assumption.

\subsection{$C^0-$Convergence Theorem}

We now recall that the well-known Gromov's compactness theorem for
cusp's curves for the compact symplectic manifolds with closed
Lagrangian submanifolds in it. For reader's convenience, we first
recall the ``weak-convergence'' for closed curves.

\vskip 3pt

{\bf Cusp-curves.} Take a system of disjoint simple closed curves
$\gamma _i$ in a closed surface $S$ for $i=1,...,k$, and denote by
$S^0$ the surface obtained from $S\setminus \cup _{i=1}^k\gamma
_i$. Denote by $\bar S$ the space obtained from $S$ by shrinking
every $\gamma _i$ to a single point and observe the obvious map
$\alpha :S^0\to \bar S$ gluing pairs of points $s'_i$ and $s''_i$
in $S^0$, such that $\bar {s}_i=\alpha (s_i')=\alpha (s_i'')\in
\bar S$ are singular (or cuspidal) points in $\bar
S$(see\cite{gro}).

\smallskip

An almost complex structure in $\bar {S}$ by definition is that in
$S^0$. A continuous map $\beta :\bar {S}\to V$ is called a
(parametrized $J-$holomorphic) cusp-curve in $V$ if the composed
map $\beta \circ \alpha :S^0\to V$ is holomorphic.

\vskip 3pt

{\bf Weak convergence.} A sequence of closed $J-$curves
$C_j\subset V$ is said to weakly converge to a cusp-curve $\bar
{C}\subset V$ if the following four conditions are satisfied

(i) all curves $C_j$ are parametrized by a fixed surface $S$ whose
almost complex structure depends on $j$, say $C_j=f_j(S)$ for some
holomorphic maps
$$f_j:(S,J_j)\to (V,J)$$

(ii) There are disjoint simple closed curves $\gamma _i\in S$,
$i=1,...,k$, such that $\bar {C}=\bar {f}(\bar {S})$ for a map
$\bar {f}:\bar {S}\to V$ which is holomorphic for some almost
complex structure $\bar {J}$ on $\bar {S}$.

(iii) The structures $J_j$ uniformly $C^\infty -$converge to $\bar
J$ on compact subsets in $S\setminus \cup _{i=1}^k\gamma _i$.

(iv) The maps $f_j$ uniformly $C^\infty -$converge to $\bar f$ on
compact subsets in $S\setminus \cup _{i=1}^k\gamma _i$. Moreover,
$f_j$ uniformly $C^0-$converge on entire $S$ to the composed map
$S\to \bar {S}\stackrel {\bar {f}}{\rightarrow }V$. Furthermore,

$$Area_{\mu }f_j(S)\to Area _{\mu }\bar {f}(\bar {S})\ for \
j\to \infty ,$$ where $\mu $ is a Riemannian metric in $V$ and
where the area is counted with the geometric
multiplicity(see\cite{gro}).

\vskip 3pt

{\bf Gromov's Compactness theorem for closed curves.} Let $C_j$ be
a sequence of closed $J-$curves of a fixed genus in a compact
manifold $(V, J, \mu ).$ If the areas of $C_j$ are uniformly
bounded,
$$Area _{\mu }\leq A,\ j=1,..,$$
then some subsequence weakly converges to a cusp-curve $\bar C$ in
$V$.

\vskip 3pt

{\bf Cusp-curves with boundary.} Let $T$ be a compact complex
manifold with boundary of dimension $1$(i.e., it has an atlas of
holomorphic charts onto open subsets of $C$ or of a closed half
plane). Its double is a compact Riemann surface $S$ with a
natureal anti-holomorphic involution $\tau $ which exchanges $T$
and $S\setminus T$ while fixing the boundary $\partial T$.
IF$f:T\to V$ is a continous map, holomorphic in the interior of
$T$, it is convenient to extend $f$ to $S$ by
$$f=f\circ \tau$$
Take a totally real submanifold $W\subset (V,J)$ and consider
compact holomorphic curves $C\subset V$ with boundaries, $(\bar C,
\partial \bar C)\subset (V,W)$, which are, topologically speaking,
obtained by shrinking to points some (short) closed loops in $C$
and also some (short) segments in $C$ between boundary points.
This is seen by looking on the double $C\cup _{\partial C}C$.

\vskip 3pt

{\bf Gromov's Compactness theorem for curves with boundary.} Let
$V$ be a closed Riemannian manifold, $W$ a totally real closed
submanifold of $V$. Let $C_j$ be a sequence of $J-$curves with
boundary in $W$ of a fixed genus in a compact manifold $(V, J, \mu
)$. If the areas of $C_j$ are uniformly bounded,
$$Area _{\mu }\leq A,\ j=1,..,$$
then some subsequence weakly converges to a cusp-curve $\bar C$ in
$V$.

\vskip 3pt

The proofs of Gromov's compactness theorem can found in
\cite{al,gro}. In our case the Lagrangian submanifold $W$ is not
compact, Gromov's compactness theorem can not be applied directly
but its proof is still effective since the $W$ has the special
geometry. In the following we modify Gromov's proof to prove the
$C^0-$compactness theorem in our case.

\vskip 3pt

  Now we state the $C^0-$convergence theorem in our case.

\begin{Theorem}
Let $(V, J, \omega , \mu )$ and $W$ as in section4. Let $C_j$ be a
sequence of $\bar J_\delta -$holomorphic section $v_j=(id,
((a^1_j,u^1_j),(a^2_j,u^2_j),f_j)):D\to D\times V$ with
$v_j:\partial D \to
\partial D\times W$ and $v_j(1)=(1,p)\in
\partial D\times W$.
constructed from section 4.  Then the areas of $C_j$ are uniformly
bounded,i.e.,
$$Area _{\mu }(C_j)\leq A,\ j=1,..,$$
and some subsequence weakly converges to a cusp-section $\bar C$
in $V$(see\cite{aul,gro}).
\end{Theorem}
Proof. We follow the proofs in \cite{gro}. Write
$v_j=(id,(a^1_{j},u^1_{j}),(a^2_j,u^2_j),f_j))$ then
$|a^2_{ij}|\leq a_0$ by the ordinary Monotone inequality of
minimal surface without boundary, see following Proposition 7.1.
Similarly $|f_j|\leq R_1$ by using the fact $f_j(\partial D)$ is
bounded in $B_{R_1}(0)$ and $\int _D|\nabla f_j|\leq 4\pi R^2$ via
monotone inequality for minimal surfaces. So, we assume that
$v_j(D)\subset V_c$ for a compact set $V_c$.

\vskip 2pt

1. {\it Removal of a net}.

\smallskip

1a.  Let $\bar V=D\times V$ and $v_j$ be regular curves. First we
study induced metrics $\mu _j$ in $v_j$. We apply the ordinary
monotone inequality for minimal surfaces without boundary to small
concentric balls $B_\varepsilon \subset (A_j,\mu _j)$ for
$0<\varepsilon \leq \varepsilon _0$ and conclude by the standard
argument to the inequality
$$Area (B_\varepsilon )\geq \varepsilon ^2, \ for \ \varepsilon
\leq \varepsilon _0;$$ Using this we easily find a interior
$\varepsilon -$net $F_j\subset (v_j,\mu _j)$ containing $N$ points
for a fixed integer $N=(\bar V, \bar J,\mu )$, such that every
topological annulus $A\subset v_j\setminus F_j$ satisfies
\begin{equation}
Diam_\mu A\leq 10length_\mu \partial A.
\end{equation}
Furthermore, let $A$ be conformally equivalent to the cylinder
$S^1\times [0,l]$ where $S^1$ is the circle of the unit length,
and let $S_t^1\subset A$ be the curve in $A$ corresponding to the
circle $S^1\times t$for $t\in [0,l]$. Then obviously
\begin{equation}
\int _a^b(lengthS_t^1)^2dt\leq Area(A)\leq C_5.
\end{equation}
for all $[a,b]\subset [0,l]$. Hence, the annulus $A_t\subset A$
between the curves $S^1_t$ and $S^1_{l-t}$ satisfies
\begin{equation}
diam_\mu A_t\leq 20 ({{C_5}\over {t}})
\end{equation}
for all $t\in [0,l]$.

\vskip 3pt

$1b$.  We consider the sets $\partial v_j\cap ((S^{1\pm})\times
F(W'\times I_i^\pm)), i=0,1$. By the construction of Gromov's
figure eight, there exists a finite components, denote it by
\begin{equation}
\partial v_j\cap ((S^{1\pm})\times F({\cal L}\times R\times I_i^\pm))
=\{\bar \gamma _{ij}^{k}\},i=0,1.
\end{equation}
we choose one point in $\bar \gamma _{ij}^{k}$ as a boundary
puncture point in $\partial v_j$ for each $i,k$.

Consider the concentric $\varepsilon $ half-disks or quadrature
$B_\varepsilon (p)$ with center $p$ on ${\bar \gamma }_{ij}^k$,
then
\begin{equation}
Area(B_\varepsilon (p))\geq \tau _0
\end{equation}
Since $Area(v_j)\leq E_0$, there exists a uniform finite puncture
points.

So, we find a boundary net $G_j\subset \partial v_j$ containing
$N_1$ points for a fixed integer $N_1(\bar V, \bar J,\mu )$, such
that every topological quadrature or half annulus $B\subset
v_j\setminus \{F_j,G_j\}$ satisfies
\begin{eqnarray}
&&\partial ''B=\partial B\cap \bar W \subset (S^{1\pm})\times
F({\cal L}\times R\times S_i), i=1,2,3,4.
\end{eqnarray}

\vskip 2pt

2. {\it Poincare's metrics}. 2a. Now, let $\mu _j^*$ be a metric
of constant curvature $-1$ in $v_j(D)\setminus F_j\cup G_j$
conformally equivalent to $\mu _j$. Then for every $\mu _j^*-$ball
$B_\rho $ in $v_j\setminus F_j\cup G_j$ of radius $\rho \leq 0.1$,
there exists an annulus $A$ contained in $v_j\setminus F_j\cup
G_j$ such that $B_\rho \subset A_t$ for $t=0.01|log|$(see Lemma
3.2.2in \cite[chVIII]{al}). This implies with $(6.3)$ the uniform
continuity of the (inclusion) maps $(v_j\setminus F_j,\mu _j^*)\to
(\bar V,\bar \mu )$, and hence a uniform bound on the $r^{th}$
order differentials for every $r=0,1,2,...$.

2b. Similarly, for every $\mu _j^*-$half ball $B_\rho ^+$ in
$v_j\setminus F_j\cup G_j$ of radius $\rho \leq 0.1$, there exists
a half annulus or quadrature $B$ contained in $v_j\setminus
F_j\cup G_j$ such that $B_\rho ^+\subset B$ with
\begin{eqnarray}
&&\partial ''B=\partial B\cap \bar W \subset (S^{1\pm})\times
F({\cal L}\times R\times S_i), i=1,2,3,4.
\end{eqnarray}
Then, by Gromov's Schwartz Lemma, i.e., Proposition 6.1-6.3
implies the uniform bound on the $r^{th}$ order differentials for
every $r=0,1,2,...$.

\vskip 2pt

3. {\it Convergence of metrics}. Next, by the standard (and
obvious ) properties of hyperbolic surfaces there is a
subsequence(see\cite{aul}), which is still denoted by $v_j$, such
that

\vskip 3pt

$(a)$. There exist $k$ closed geodesics or geodesic arcs with
boundaries in $\partial v_j\setminus F_j$, say
$$\gamma _i^j\subset (v_j\setminus F_j,\mu ^*_j), i=1,...,k,j=1,2,...,$$
whose $\mu _j^*-$length converges to zero as $j\to \infty $, where
$k$ is a fixed integer.

\vskip 2pt

$(b)$. There exist $k$ closed curves or geodesic arcs with
boundaries in $\partial S$ of  a fixed surface, say $\gamma _j$ in
$S$, and an almost complex structure $\bar J$ on the corresponding
(singular)  surface $\bar S$, such that the almost complex
structure $J_j$ on $v_j\setminus F_j$ induced from $(V,J)$
$C^\infty -$converge to $\bar J$ outside $\cup _{j=1}^k\gamma _j$.
Namely, there exist continuous maps $g^j:v_j\to \bar S$ which are
homeomorphisms outside the geodesics $\gamma _i^j$, which pinch
these geodesics to the corresponding singular points of $\bar
S$(that are the images of $\gamma _i$) and which send $F_j$ to a
fixed subset $F$ in the nonsingular locus of $\bar S$. Now, the
convergence $J_j\to \bar J$ is understood as the uniform $C^\infty
-$convergence $g^j_*(J_j)\to \bar J$ on the compact subsets in the
non-singular  locus $\bar S^*$ of $\bar S$ which is identified
with $S\setminus \cup _{i=1}^k\gamma _i$.

4. {\it $C^0-$interior convergence}.  The limit cusp-curve $\bar
v:\bar S ^*\to  \bar V$, that is a holomorphic map which is
constructed by first taking the maps
$$\bar v_j=(g_j)^{-1}:S\setminus \cup _{i=1}^k\gamma _i\to \bar V$$
Near the nodes of $\bar S$ including interior nodes and boundary
nodes, by the properties of hyperbolic metric $\mu ^*$ on $\bar
S$, the neighbourhoods of interior nodes are corresponding to the
annulis of the geodesic cycles. By the reparametrization of $v_j$,
called $\bar v_j$ which is defined on $S$ and extends the maps
$\bar v_j:S\to S_j\to V$(see\cite{aul,gro}). Now let $\{
z_i|i=1,...,n\}$ be the interior nodes of $\bar S$. Then the
arguments in \cite{aul,gro} yield the $C^0-$interir convergece
near $z_i$.

5. {\it $C^0-$boundary convergence}. Now it is possible that the
boundary of the cusp curve $\bar v$ does not remain in $\bar W$.
Write $\bar v(z)=(h,((a_1,u_1),(a_2,u_2)),f)(z)$, here $h(z)=z$ or
$h(z)\equiv z_i$, $i=1,...,n$, $z_i$ is cusp-point or bubble
point. We can assume that $\bar p=(1,p)\in \bar v_n$ is a puncture
boundary point. Let $\bar v_1$ be the component of $\bar v$ which
through the point $\bar p$. Let $D=\{z|z=re^{i\theta},0\leq
r\leq1,0\leq \theta \leq 2\pi\}$. We assume that $\bar
v_1:D\setminus \{e^{i\theta _i}\}_{i=1}^k\to V_c$, here
$e^{i\theta _i}$ is node or puncture point. Near $e^{i\theta _i}$,
we take a small disk $D_i$ in $D$ containing only one puncture or
node point $e^{i\theta _i}$. By the reparametrization and the
convergence procedure, we can assume that $\bar v_{1i}=(\bar
v_1|D_i)$ as a map from $D^+\setminus \{0\}\to V_c$ with $\bar
v_1([-1,1]\setminus \{0\})\subset S^1\times F(W'\times S^1)$ and
$area (\bar v_{1i})\leq a_0$, $a_0$ small enough.
 Since $Area(\bar v_{1i})\leq a_0$, there exist curves $c_k$
near $0$ such that $l(\bar v_{1i}(c_k))\leq \delta _1$. By the
construction of convergence, we can assume that $l(\bar
v_n(c_k))\leq 2\delta _1$. If $\bar v_{1i}(\partial c_k)\subset
(S^{1})\times F({\cal L}\times [-N_0,N_0]\times S^1)$, we have
$\bar v_n(\partial c_k)\subset (S^{1})\times F({\cal L}\times
[-2N_0,2N_0]\times S^1)$ for $n$ large enough. Now $\bar v_n(c_k)$
cuts $\bar v_n(D)$ as two parts, one part corresponds to $\bar
v_{1i}$, say $\bar u_n(D)$. Then $area (\bar u_n(D))=
area(h_{n1})+|\Psi '(u_{n2}(c_k^1))-\Psi '(u_{n2}(c_k^2))|$, here
$\partial c_k=\{c_k^1,c_k^2\}$. Then by the proof of Lemma6.1-6.3,
we know that $\bar u_n(\partial D\setminus c_k)\subset
(S^{1})\times F({\cal L}\times [-100N_0,100N_0]\times S^1)$. So,
$\bar v_{1i}([-1,1]\setminus \{0\})\subset S^1\times F({\cal
L}\times [-100N_0,100N_0]\times S^1)$. By proposition 6.6, one
singularity of $\bar v_1$ is deleted. We repeat this procedure, we
proved that $\bar v_1$ is extended to whole $D$. So, the boundary
node or puncture points of $\bar v$ are removed. Then by choosing
the sub-sub-sequences of $\mu ^*_j$ and $\bar v_j$, we know that
$\bar v_j$ converges to $\bar v$ in $C^0$ near the boundary node
or puncture point. This proved the $C^0-$boundary convergence.
Since $\bar v_j(1)=\bar p$, $\bar p\in \bar v(\partial D)$, $\bar
v(\partial D)\subset \bar W$.

\vskip 2pt

6. {\it Convergence of area}. Finally by the $C^0-$convergence and
$area(v_j)=\int _Dv_j^* \bar \omega $, one easily deduces
$$area(v(S))=\lim _{j\to \infty }(v_j(S_j)).$$

\subsection{Bounded image of $J-$holomorphic curves in $W$}

\begin{Proposition}
Let $v$ be the solutions of equations (4.16), then
$$d_W(p,v(\partial D^2))
=max\{ d_W(p,q)|q\in f(\partial D^2)\} \leq d_0<+\infty$$
\end{Proposition}
Proof. It follows directly from
Gromov's $C^0-$convergence theorem.

\section{Proof of Theorem 1.1}

\begin{Proposition}
If $J-$holomorphic curves $C\subset \bar V$ with
boundary
$$\partial C\subset D^2\times ([0,\varepsilon ]\times \Sigma )
\times ([0,\varepsilon ]\times \Sigma )\times R^2$$
and
$$C\cap (D^2\times (\{-3\}\times \Sigma )
\times (R\times \Sigma )\times R^2)\ne \emptyset $$
or
$$C\cap (D^2\times (R\times \Sigma )
\times (\{-3\}\times \Sigma )\times R^2)\ne \emptyset $$
Then
$$area(C)\geq 2l_0.$$
\end{Proposition}
Proof. It is obvious by monotone inequality argument for
minimal surfaces.

\begin{Note}
we first observe that
any $J-$holomorphic curves with boundary in
$R^+\times \Sigma $ meet the
hypersurface $\{-3\}\times \Sigma $ has
energy at least $2l_0$, so we take
$\varepsilon $ small enough such that
the Gromov's figure eight contained
in $B_{R(\varepsilon )}\subset C$ for $\varepsilon $ small enough and
the energy
of solutions in section 4 is smaller than
$l_0$. we specify the constant
$a_0$, $\varepsilon $ in section 3.1-3 such that
the above conditions satisfied.
\end{Note}

\begin{Theorem}
There exists a non-constant $J-$holomorphic map $u: (D,\partial D)\to
(V'\times C,W)$ with $E(u)\leq 4\pi R(\varepsilon )^2$
for $\varepsilon $ small enough such that $4\pi R(\varepsilon )^2\leq
l_0$.
\end{Theorem}
Proof.  By Proposition 5.1, we know that the image $\bar v(D)$
of solutions of equations (\ref{eq:4.16}) remains a bounded or
compact part of the non-compact Lagrangian submanifold
$W$. Then, all arguments in \cite{aul,gro} for the case $W$ is closed
in $S\Sigma \times S\Sigma \times R^2$ can be extended to our case, especially
Gromov's $C^0-$converngence theorem holds. But the results in
section 4 shows the solutions of equations (\ref{eq:4.16}) must
denegerate to a cusp curves, i.e., we obtain a Sacks-Uhlenbeck-Gromov's
bubble, i.e., $J-$holomorphic sphere or disk with boundary
in $W$, the exactness of $\omega $ rules out
the possibility of $J-$holomorphic sphere. For the more detail, see the proof
of Theorem 2.3.B in \cite{gro}.

\vskip 3pt

{\bf Proof of Theorem 1.1}.
If $(\Sigma ,\lambda )$ has not closed Reeb orbit, then
we can construct a Lagrangian
submanifold $W$ in $V=V'\times C$, see section 3.
Then as in section 4, we construct an
anti-holomorphic section $c$ and
for large vector $c\in C$ we know
that the nonlinear
Fredholm section or Cauchy-Riemann section
has no solution, this implies that the section is non-proper, see section 4.
The non-properness of the
section and
the Gromov's compactness
theorem in section 6 implies the existences of the cusp-curves.
So, we must have the $J-$holomorphic sphere or
$J-$holomorphic disk with bounadry in $W$.
Since the symplectic manifold $V$ is
an exact symplectic mainifold and $W$ is an exact
Lagrangian submanifold in $V$, by Stokes formula,
we know that
the possibility of $J-$holomorphic sphere or disk
elimitated.
So our priori assumption does not hold which
implies the contact maifold $(\Sigma ,\lambda )$
has at least closed Reeb orbit. This finishes
the proof of Theorem 1.1.


\begin{thebibliography}{99}



\bibitem{ar} Arnold, V. I., First steps in symplectic topology, Russian Math.
Surveys 41(1986),1-21.

\bibitem{arg}  Arnold, V.\& Givental, A., Symplectic Geometry,
in: Dynamical Systems IV, edited by V. I. Arnold and S. P. Novikov,
Springer-Verlag, 1985.

\bibitem{aul}  Audin, M\& Lafontaine, J., eds.: Holomorphic Curves in Symplectic
Geometry. Progr. Math. 117, (1994) Birkha\"{u}ser, Boston.

\bibitem{ch}  Chaperon, M., Questions de geometrie symplectique,
in Seminaire Bourbaki, Asterisque 105-106(1983), 231-249.


\bibitem{el} Eliashberg, Y., Symplectic topology in the nineties,
Differential geometry and its applications 9(1998)59-88.

\bibitem{elg} Eliashberg,Y.\& Gromov, M., Lagrangian Intersection
Theory: Finite-Dimensional Approach, Amer. Math. Soc. Transl. 186(1998):
27-118.

\bibitem{em} Eliashberg,Y.\& Mishachev, N.M., Holonomic approximation
and Gromov's h-principle, 15pages,
AIM2001-3, SG(GI).


\bibitem{ehs} Eliashberg,Y., Hofer,H., \& Salamon,S., Lagrangian Intersections
in contact gemetry, Geom. and Funct. Anal., 5(1995): 244-269.

\bibitem{fl} Floer, A., Symplectic fixed point and holomorphic
spheres, Commun. Math. Phys., 120(1989), 575-611

\bibitem{fhv} Floer, A., Hofer, H.\& Viterbo, C., The Weinstein
conjecture in $P\times C^l$, Math.Z. 203(1990)469-482.


\bibitem{gi} Givental, A. B., Nonlinear generalization of the
Maslov index, Adv. in Sov. Math., V.1,
AMS, Providence, RI, 1990.

\bibitem{gra} Gray, J.W., Some global properties of contact
structures. Ann. of Math., 2(69): 421-450, 1959.

\bibitem{gro}  Gromov, M., Pseudoholomorphic Curves in Symplectic manifolds.
Inv. Math. 82(1985), 307-347.

\bibitem{grob}  Gromov, M., Partial
Differential Relations, Springer, New York, 1986.


\bibitem{ho}  Hofer, H., Pseudoholomorphic curves in symplectizations with
applications to the Weinstein conjection in dimension three.
Inventions Math., 114(1993), 515-563.


\bibitem{hv} Hofer, H.\& Viterbo, C., The Weinstein
conjecture in the presence of holomorphic spheres,
Comm. Pure Appl. Math.45(1992)583-622.

\bibitem{hz} Hofer, H.\& Zehnder, E., Periodic solutions on hypersurfaces and
a result by C. Viterbo, Invent. Math. 90(1987)1-9.



\bibitem{hor} H\"ormander, L., The Analysis of Linear Partial
Differential Operators I, Springer-Verlag, 1983.


\bibitem{hu} Hummel, C., Gromov's Compactness Theorem for
Pseudo-holomorphic Curves, Progress In Math. Volume 151, Birkhauser Verlag, 1997.



\bibitem{kl}  Klingenberg, K., Lectures on closed Geodesics, Grundlehren der Math.
Wissenschaften, vol 230, Spinger-Verlag, 1978.


\bibitem{ls} Lalonde, F \& Sikorav, J.C., Sous-Vari\`t\`es Lagrangiennes
et lagrangiennes exactes des fibr\`es cotangents, Comment. Math. Helvetici 66(1991)
18-33.


\bibitem{ma1}  Ma, R., A remark on the Weinstein conjecture in $M\times R^{2n}$.
Nonlinear Analysis and Microlocal Analysis, edited by K. C. Chang,
Y. M. Huang and T. T. Li, World Scientific Publishing, 176-184.

\bibitem{ma2}  Ma, R., Symplectic Capacity and Weinstein Conjecture in
Certain Cotangent bundles and Stein manifolds. NoDEA.2(1995):341-356.

\bibitem{ma3} Ma, R., Legendrian submanifolds and A Proof on
Chord Conjecture, Boundary Value Problems, Integral
Equations and Related Problems, edited by J K Lu \& G C Wen,
World Scientific, 135-142.

\bibitem{ma4} Ma, R., The existence of
$J-$holomorphic curves and applications to the Weinstein
conjecture, Chin. Ann. of Math. 20B:4 (1999), 425-434.

\bibitem{ra} Rabinowitz, P., Periodic solutions of Hamiltonian systems,
Comm. Pure. Appl. Math 31, 157-184, 1978.


\bibitem{sau} Sacks, J. and Uhlenbeck,K., The existence of
minimal 2-spheres. Ann. Math., 113:1-24, 1983.

\bibitem{sm} Smale, S., An infinite dimensional
version of Sard's theorem, Amer. J. Math. 87: 861-866, 1965.

\bibitem{thu} Thurston, B., The theory of foliations in codimension greater than
one, Comm. Math. Helv. 49: 214-231, 1974.

\bibitem{vi}  Viterbo, C., A proof of the Weinstein conjecture in
$R^{2n} $, Ann. Inst. Henri. Poincar$\acute e$, Analyse
nonlin$\grave e$aire, 4: 337-357, 1987.

\bibitem{vi1} Viterbo, C., Exact Lagrange submanifolds, Periodic
orbits and the cohomology of free loop spaces, J.Diff.Geom.,
47(1997), 420-468.

\bibitem{we1}  Weinstein, A., Periodic orbits for convex
Hamiltonian systems. Ann. Math. 108(1978),507-518.


\bibitem{we2}  Weinstein, A., On the hypothesis of Rabinowitz's periodic
orbit theorems, J. Diff. Eq.33, 353-358, 1979.


\bibitem{wen}  Wendland, W., Elliptic systems in the plane, Monographs and
studies in Mathematics 3, Pitman, London-San Francisco, 1979.

\end{thebibliography}
\end{document}